\documentclass[12pt]{article}
\usepackage[utf8]{inputenc}
\usepackage{fullpage}
\usepackage{amsmath}
\usepackage{amsfonts}
\usepackage{amssymb}
\usepackage{graphicx}
\usepackage{subcaption}
\usepackage{xcolor}
\usepackage{url}

\newcommand{\R}{\mathbb{R}}

\newcommand{\cA}{\mathcal{A}}

\usepackage{todonotes}

\graphicspath{{Figures/}}

\title{Evaluating the skill of a geometric early warning for tipping in a rapidly forced nonlinear system}

\author{Paul Ritchie\thanks{Department of Mathematics and Statistics, University of Exeter, Exeter EX4 4QF, UK}\and Sneha Kachhara${}^*$ \and Peter Ashwin${}^*$}

\date{\today}

\begin{document}

\maketitle

\begin{abstract}
The future behavioural fate of a forced nonlinear system can depend sensitively on the forcing profile as well as natural fluctuations within the system. This is especially the case for rate-induced tipping, where the forcing pushes the system to a basin boundary of a future behaviour and small changes in the forcing can lead to drastically different eventual behaviours. This sensitivity may be present only for a limited period of time, for example when the forcing is most rapidly changing. Moreover, critical slowing down based methods fail to be informative in such cases. We investigate a geometric early warning to evaluate when a system is in such a sensitive state. This involves computing the R-tipping indicator, namely the signed distance to an approximate R-tipping threshold. The latter is a dynamic state that embeds knowledge of the system and future behaviour of the forcing.  We contrast this with early warnings of bifurcation-induced tipping, where tipping is associated with passing a threshold on slow variation of forcing. As an example, we consider methods of early prediction of the future state for a 3-box model of the Atlantic Meridional Overturning Circulation (AMOC) with specified rapid forcing. We show that the skill of the geometric early warning compares favourably with simple thresholds.
\end{abstract}

\tableofcontents

\section{Introduction}

In many application areas (e.g. for climate \cite{armstrong2022exceeding,Lenton08,clarke2021compost}, ecology \cite{dakos2019ecosystem,feudel2023rate,scheffer2008pulse}, neuroscience \cite{trevelyan2025brain}, and finance \cite{smug2022adaptive}) dynamical models need to be not only nonlinear (meaning that they can show behaviours such as multistability and chaos), but also nonautonomous (meaning that the dynamic rules governing the system evolution depend on some time-varying input - also called the forcing). If the time-varying input is slow relative to the system timescale then one can understand much of the behaviour in terms of quasi-static (i.e. asymptotically slow relative to system timescale) variation of attractors of a {\em frozen} autonomous system depending on an input parameter. Passing through a bifurcation will cause an abrupt response in the system if an attractor that is being followed quasi-statically loses stability or disappears at a bifurcation point \cite{ashwin2012tipping,ashwin2017parameter,wieczorek2023}.

Early warning indicators have been developed to forewarn the approach of such B-tipping events \cite{dakos2024tipping,bury2021deep,lenton2011early}. As a bifurcation is approached, the stability of the current state will weaken, meaning that fluctuations become more pronounced and the system will take longer to recover back to the stable state, a phenomenon known as \textit{critical slowing down} \cite{lenton2011early,scheffer2012anticipating}. Various early warning indicators, such as observing an increase in lag-1 autocorrelation and variance, have been designed to detect critical slowing down, and to provide forewarning of the approaching bifurcation or tipping point \cite{dakos2008slowing,dakos2024tipping}. For forcing that overshoots a bifurcation threshold, there is an asymptotic inverse square law that determines by how much and for how long a bifurcation threshold can be exceeded before tipping will be triggered \cite{ritchie2019inverse}, although applying this theory as a warning is not so easy. 

If the forcing is rapid (relative to system timescale), standard early warning indicators based on critical slowing down  will fail \cite{van2021no,clarke2026conds,chapman2025quantifying} or at least there will be a delay that means the warning is no longer early \cite{ritchie16}. This is because a rapidly forced system can explore states far from any attractor of the quasi-static or frozen system. This makes detection of declining resilience of the attractor more challenging or even impossible. Indeed, it is known that a rapidly forced system can avoid tipping even if the tipping point or fold bifurcation is crossed temporarily \cite{ritchie21}, or may undergo rate-induced tipping (R-tipping) even if there is no bifurcation present for slow forcing in the same range \cite{ashwin2012tipping,wieczorek2023}. 

Some attempts have been made to develop early warnings in the presence of rapid forcing that overcome these difficulties. In \cite{huang2024anticipating}, a deep learning approach is taken to identify the members of the ensemble that tip, and this can have predictive skill, although it can be hard to explain precisely which features are important to achieve this skill. A recent paper \cite{chapman2025quantifying} proposed taking a {\em geometric early warning} approach to understand the future fate of systems that have been through a period of rapid transient forcing. 

The approach in \cite{chapman2025quantifying} aims to find a threshold in phase space that corresponds to the R-tipping threshold \cite{wieczorek2023}. This threshold corresponds to the set of trajectories that separate different attractors in the future limit; this is a notion of attractor basin boundary adapted to a setting that is nonautonomous (but with an autonomous future limit). The R-tipping threshold is a boundary in the phase space that separates tipping from not tipping based on the remainder of the forcing profile; however this threshold is dynamic - it moves with forcing. Further complication is introduced if the system considered is stochastic, in which case the rate- and noise- dependent aspects of the forcing will act together \cite{ritchie17,slyman2023rate}, or if there is a mixture of rate-dependent and chaotic forcing \cite{ashwin2025contrasting}.

In this paper, we extend the ideas in \cite{chapman2025quantifying} and propose a geometric early warning signal based on measuring the signed distance to the R-tipping threshold, which works even when the forcing is still underway. This requires knowledge of the future forcing, system dynamics, and initial condition. This has the advantage of being fully explainable (though still subject to uncertainties in the system dynamics, future forcing, and initial state). 

Suppose a dynamical system is subjected to an external forcing profile (or an ensemble of profiles with quantified probability density) for a known initial state. Knowing where the system is relative to its R-tipping threshold will help indicate whether the system is committed to one future behaviour or another.  In the presence of noise, we need probabilistic estimates that depend on the initial state in phase space, the underlying system dynamics, and the remaining profile of the forcing. However, note that we do not need to know anything about the past (such as the past forcing profile and/or how the system reached its current position), although given knowledge of the past forcing, we can expect the system to be close to a pullback attractor which may be a single trajectory of the time-dependent system \cite{ashwin2017parameter}, or a nontrivial time-dependent set that is a union of trajectories \cite{alkhayuon2018rate,lohmann2024predictability}.

The paper is structured as follows: after presenting a motivating example in Section~\ref{sec:mot_example}, in Section~\ref{sec:Rtipping_edge_states} we recall the precise definition of the R-tipping threshold. We discuss how it can be approximated and how the signed distance to this threshold can be used as an early warning signal. In Section~\ref{sec:AMOC_example} we use a strongly forced conceptual box model for the AMOC as an illustrative example. We compare the R-tipping indicator with simpler thresholds in phase space as early warning signals for rate-induced tipping. Quantification of the skill of these early warning signals is considered in Section~\ref{sec:EWS_skill}. We discuss extensions and potential applications of these results in Section~\ref{sec:Discussion}.   

\subsection{A motivating example: temporary loss of an attractor}
\label{sec:mot_example}

Consider a (possibly stochastic) system with two attractors such that there is a temporary loss of one of these attractors. This can be thought of as a brief overshoot of a fold \cite{ritchie2019inverse,ritchie2023rate}, but where the forcing is monotonically increasing over time. The system for $X(T)\in\R$ we consider is
\begin{equation}
\mathrm{d}x = f(x,P_{\theta}(t))\mathrm{d}t+\sigma\mathrm{d}W_t,
\label{eq:generic_ex}
\end{equation}
where
$$
f(x,p):= -\frac{x^3}{3} + x - p(p_+-p).
$$
We assume additive white noise of strength $\sigma\geq 0$, expressed in terms of a Wiener process $W_t$. We assume the external forcing is a monotonic parameter shift
\begin{equation}
P_{\theta}(t) = \frac{p_+}{2}\left(\tanh(\theta t)+1\right),
\label{eq:ramp_forcing}
\end{equation}
which starts at $0$ and asymptotically approaches some maximal level $p_+>0$. The rate of external forcing can be scaled using the parameter $\theta\geq 0$. Note that the frozen deterministic system (corresponding to $\sigma=0$ and $\theta\rightarrow 0$) is simply
\begin{equation}
\dot{x} = f(x,p)
\label{eq:frozen_ex}
\end{equation}
for $0\leq p \leq p_+$. For most values of $p$ there are two attractors that we term the upper and lower attractors. There is a critical value $p_c=\frac{2\sqrt{2}}{\sqrt{3}}=1.6330$ such that if $p_+>p_c$, the upper attractor disappears for some interval in the interior of  $(0,p_+)$ through a pair of fold bifurcations. We assume we start on the upper attractor and are concerned whether the state moves to the lower attractor, in which case we say it {\em tips}. 

\begin{figure}[htp!]
\centering
\includegraphics[width=\linewidth]{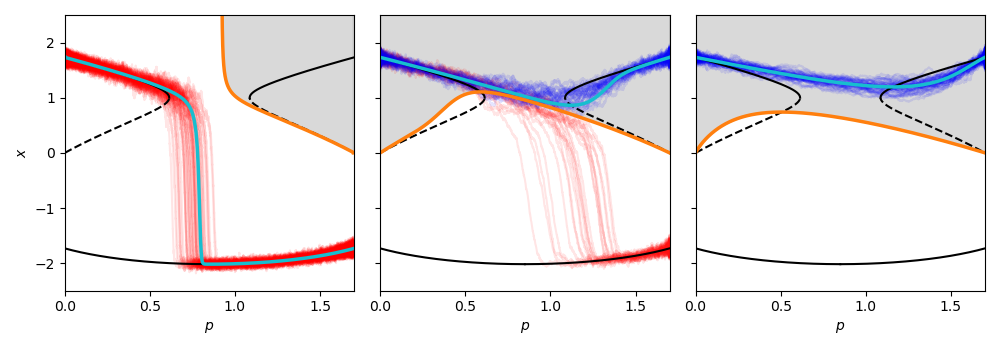}%
\caption{\textbf{Rate and noise-induced effects for temporary loss of an attractor.} 
Simulations of \eqref{eq:generic_ex}--\eqref{eq:ramp_forcing} showing (left) slow ($\theta=0.015$) forcing, (middle) intermediate speed ($\theta=0.08$), (right) fast ($\theta = 0.4$).
Stable and unstable equilibria of the frozen system are given by the black solid and dashed curves respectively and we set $p_+=1.7>p_c$ so that there is a finite range of times where there is only one attractor for the frozen system.
The deterministic case ($\sigma=0$) shows the local pullback attractor as a light blue curve. The R-tipping threshold (shown as the orange curves) divides the plane into points in phase space that undergo tipping to the lower attractor (white region) from those that continue to the upper attractor (grey region). In the noisy case ($\sigma=0.1$) some of an ensemble of trajectories (red) tip to the lower attractor but the others (blue) continue to the upper attractor. In the intermediate noisy case there is a splitting of the ensemble.}
\label{fig:slow_fast_overshoot_Rtipping_threshold_traj}
\end{figure}

In the noise-free/deterministic case ($\sigma=0$), if $p_+>p_c$ then the value of $\theta$ will determine whether there is tipping onto the lower attractor; see Figure~\ref{fig:slow_fast_overshoot_Rtipping_threshold_traj}.  For sufficiently slow ($\theta$ small) change in external forcing, the local pullback attractor (i.e. the deterministic ($\sigma=0$) solution of \eqref{eq:generic_ex}--\eqref{eq:ramp_forcing} initialised on the upper attractor: see \cite{ashwin2017parameter}) will track the upper stable branch of equilibria of the frozen system up to the first fold bifurcation. At the fold, this upper attractor is lost, causing the pullback attractor to tip to the lower stable branch before the upper attractor reappears at higher values of $p$; see left panel of Figure~\ref{fig:slow_fast_overshoot_Rtipping_threshold_traj}. The orange curve shows the R-tipping threshold \cite{wieczorek2023}. This divides those initial conditions that go to the upper attractor in the future (grey region) and those that go to the lower attractor (white region).
If the change in forcing is sufficiently fast ($\theta$ large), then the system will avoid tipping to the lower attractor: the local pullback attractor (see the right panel of Figure~\ref{fig:slow_fast_overshoot_Rtipping_threshold_traj}). The middle panel is close to a critical value of $\theta=\theta_c$ depending on the choice of $p_+$ such that the pullback attractor and the R-tipping threshold coincide. This scenario can be related to safe overshoot \cite{ritchie21} and rate-induced tracking \cite{duenas2023rate} of a non-monotonically forced system. 

In the noisy case ($\sigma>0$) the deterministic picture becomes blurred, especially near the critical rate $\theta_c$. Figure~\ref{fig:slow_fast_overshoot_Rtipping_threshold_traj} show members of an ensemble of realizations that eventually find the lower state (red trajectories) or the upper state (blue trajectories). Near the critical rate $\theta\approx \theta_c$ for the deterministic case, observe that the ensemble may split, as shown in the middle panel.
This example highlights that leaving the basin of attraction of the original state is not sufficient to signal tipping. Indeed, certain forcing scenarios will result in reduced predictability of long-term behaviour \cite{lohmann2024predictability,lohmann2025role} due to the ensemble of possible system trajectories splitting into more than one storyline. This scenario is shown in the middle panel of Figure~\ref{fig:slow_fast_overshoot_Rtipping_threshold_traj}, where the ensemble splits along an R-tipping threshold (orange curve) \cite{wieczorek2023} that corresponds to the trajectory that limits to the edge state in future time. The pullback attractor limits to the original state but comes very close to this R-tipping threshold, signifying a close to critical rate of change that separates tipping from not tipping deterministically \cite{ritchie2023rate}. However, noise has the capability to kick the system into tipping (red trajectories), whereas on other occasions tipping is avoided (blue trajectories). 

Understanding this behaviour, by being able to distinguish between tipping and non-tipping trajectories, has been an area of recent focus to anticipate R-tipping. For an R-tipping collapse of the Atlantic Meridional Overturning Circulation (AMOC), \cite{chapman2025quantifying} 
suggested the signed distance to the R-tipping threshold as a potential geometric early warning signal for tipping. We call this the {\em R-tipping indicator}.
In the study \cite{chapman2025quantifying}, the early warning was only examined after the end of transient forcing. In this paper, we demonstrate that this limitation can be removed and verify that this approach can be more skilful than others.

\section{R-tipping edge states and thresholds}
\label{sec:Rtipping_edge_states}

Consider a general nonautonomous evolution that depends on time-dependent forcing, possibly in the presence of noise. We express this as an It\^{o} stochastic differential equation for $x\in \R^d$:
\begin{equation}
dx=f(x,P_{\theta}(t))dt+\sigma dW_t
\label{eq:SDE-nonaut}
\end{equation}
where $\sigma$ is such that $\sigma^T\sigma$ is positive semidefinite (but small) and $W_t$ is a $d$-dimensional Wiener process, We assume $f:\R^d\times\R^l\rightarrow \R^d$ is smooth and call $P_{\theta}(t)\in\R^l$ the {\em forcing}. This is a parameter path in $t$ that itself is parametrized by $\theta\in\R$. 
The noise-free case ($\sigma\equiv 0$), (\ref{eq:SDE-nonaut}) reduces to a nonautonomous ordinary differential equation for $x\in\R^d$ of the form
\begin{equation}
\dot{x}=f(x,P_{\theta}(t))
\label{eq:ODE-nonaut}
\end{equation}
We call
\begin{equation}
\dot{x}=f(x,p)
\label{eq:ODE-frozen}
\end{equation}
the {\em frozen system} associated with (\ref{eq:ODE-nonaut}). R-tipping may appear in systems such as (\ref{eq:ODE-nonaut}) where there is asymptotically stationary forcing (e.g. \cite{duenas2023rate}). A special case of this is when the forcing is asymptotically constant, in which case there are autonomous limits and, as long as the forcing decays fast enough, we can compactify time to create an autonomous system on a higher-dimensional state space \cite{wieczorek2023}: we focus on this case. 

\subsection{R-tipping for nonautonomous systems with autonomous limits}

In the case that $\lim_{t\rightarrow \pm\infty} P_{\theta}(t)=p_{\pm}$ then we say $P_{\theta}$ is {\em asymptotically constant}. For simplicity of presentation, we assume $p_{\pm}$ are independent of $\theta$ and that there is only one parameter, so $l=1$. 
It is instructive to relate the dynamics of the nonautonomous system (\ref{eq:ODE-nonaut}) to the dynamics of the asymptotic systems when $p = p_{\pm}$ and, more generally, to the dynamics of the frozen systems (\ref{eq:ODE-frozen}) for other values of $p$.

\begin{figure}
\begin{center}
\includegraphics[width=15cm]{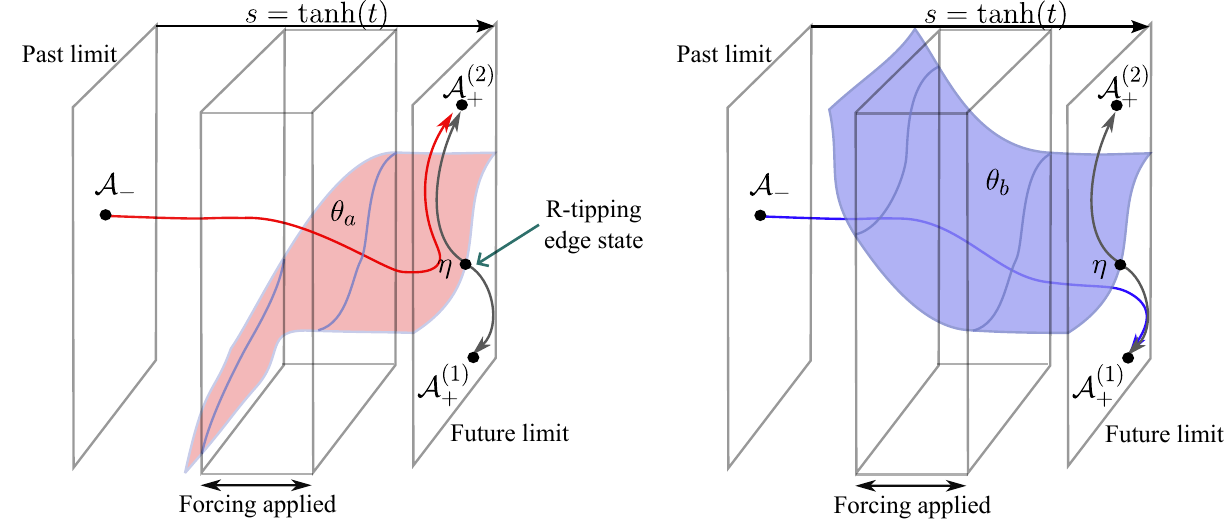}
\end{center}
\caption{\textbf{Schematic illustration of R-tipping threshold for a system with time-dependent forcing that is asymptotically constant in the past and future.} The horizontal axis represents a compactified time (say $s=\tanh(t)$) while the vertical plane represents the phase space of the system. Suppose the system starts near an attractor $\mathcal{A}_-$ in the past and has forcing applied in the boxed region with the forcing parameterized by some $\lambda$, and suppose that there is an R-tipping edge state $\eta$ that divides the basins of the attractors $\mathcal{A}_+^{(1,2)}$. The R-tipping threshold (stable manifold of $\eta$ for the compactified extended system)  is shown in red and blue for values of the parameter $\theta$ near an R-tipping in the interval $\theta_a<\theta<\theta_b$. The associated system trajectories (pullback attractors) starting at $\mathcal{A}_-$ are shown in red/blue. For $\theta$ at the point of R-tipping, the system state evolves towards the edge state $\eta$. }
\label{fig:RtippingEdge}
\end{figure}

Suppose that we have an attractor $\cA_-$ for the past limit system when $p=p_{-}$, then one can show that there is a local pullback attractor for the nonautonomous system that limits to $\cA_-$ in the past \cite{alkhayuon2018rate,ashwin2017parameter}. This pullback attractor may split into the basin of more than one attractor $\cA_+^{(k)}$ enumerated by $k$ in the future \cite{lohmann2024predictability,ashwin2021physical}. If $\cA_{-}$ is a stable equilibrium, there is a single pullback attracting trajectory, and this typically limits to one attractor in the future. However, if a parameter $\theta$ in the forcing is changed, the future limit of a pullback attractor may also change: we term this a rate-induced transition (R-tipping) even when $\theta$ is not necessarily a rate. 
For asymptotic autonomous limits, we can associate R-tipping with the presence of connecting orbits in a system, for a compact and autonomous extension of the original system. These connections will be from $\cA_-$ in the past limit to an {\em R-tipping edge state} $\eta$ in the future limit, on the boundary between two attractor basins. To do this, we define a {\em compactified time} $s$
\begin{equation}
    s=\tanh(\alpha t)
    \label{eq:compactify-t}
\end{equation} 
where $\alpha>0$ is a compactification parameter that can, under suitable assumptions, be chosen such that a compactified extended system is as smooth as the original system at the past and future limits \cite{wieczorek2023}. In this case, invariant manifold theory can be used to understand (and compute) the connecting orbits associated with R-tipping. If we write the inverse of (\ref{eq:compactify-t}) as $t=h_{\alpha}(s)$ then (\ref{eq:ODE-nonaut}) can be written in the form of a {\em compactified extended system}
\begin{equation}
    \begin{aligned}
        \dot{x}=&f(x,P_{\theta}(h_{\alpha}(s)))\\
        \dot{s}=&\alpha (1-s^2)
    \label{eq:ODE-compactified}
    \end{aligned}
\end{equation}
with $x\in\R^d$, $s\in[-1,1]$ and we take, as in \cite{wieczorek2023}, the continuous extension of $P_{\theta}(h_{\alpha}(s))$ to the limit $s=\pm 1$. We show in Figure~\ref{fig:RtippingEdge} a schematic illustration of an R-tipping threshold for an R-tipping edge state $\eta$ in the future limit as a parameter $\theta$ changes through R-tipping. 

For (\ref{eq:ODE-compactified}), note that the past limit system, which contains $\cA_-$ corresponds to the invariant set $s=-1$ while the future limit system, which contains $\cA_+^{(k)}$ corresponds to the invariant set $s=+1$. If there is more than one attractor for the future limit, then the basin boundary can contain compact invariant sets whose stable manifolds form part (or all) of the boundary. The invariant sets on the basin boundaries that have codimension one stable manifolds are edge states for the future limit system. Moreover, they will have codimension-one stable manifolds in the full compactified extended system (\ref{eq:ODE-compactified}) due to the additional contracting ($s$) direction near the subspace $s=1$. We call this stable manifold in the full system the {\em R-tipping threshold} \cite{wieczorek2023}. This R-tipping threshold is therefore an object in phase space that can be approximated by integrating backwards in time from the stable manifold of the edge state at the end of the forcing. Details of its calculation and that of the R-tipping indicator are provided in Appendix~\ref{app:R-threshold_calc}.

We can similarly consider the compactified extended version of (\ref{eq:SDE-nonaut}) as
\begin{equation}
    \begin{aligned}
        dx=&f(x,P_{\theta}(h_{\alpha}(s)))dt+ \sigma dW_w\\
        ds=&\alpha (1-s^2)dt.
    \label{eq:SDE-compactified}
    \end{aligned}
\end{equation}
and note that over timescales that are small compared to the timescale of noise-induced transitions between attractors (i.e. for sufficiently small amplitude noise) we can characterise R-tipping transitions depending on a forcing parameter $\theta$.

\section{Example: a low-order model of the AMOC}
\label{sec:AMOC_example}

As an illustration of our approach, we will use the 3-box simplification discussed in \cite{Alkhayoun2019,chapman2025quantifying} of the data-adapted 5-box model for the Atlantic Meridional Overturning Circulation (AMOC) introduced by Wood et al \cite{wood2019observable}. 
The model is studied in \cite{chapman2025quantifying} for geometric early warnings after transient forcing.

\subsection{Rapid forcing of a three-box model of AMOC}

The model we consider is for mean salinities of the global oceans in the North Atlantic ($S_{\rm{N}}$), the Tropical Atlantic ($S_{\rm{T}}$), the Southern Ocean ($S_{\rm{S}}$), the Bottom Ocean ($S_B$) and the Indo-Pacific ($S_{\rm{IP}}$). A reduced version of this model eliminates the salinity evolution of two of the ocean boxes, $S_B$ and $S_S$, by assuming that they are in equilibrium. The salinity in the Indo-Pacific box $S_{\rm{IP}}$ can be further eliminated by assuming the total salinity content $C$ is conserved:
\begin{equation}
	S_{\rm{IP}} = \{C- (V_{\rm{N}} S_{\rm{N}}+ V_{\rm{T}} S_{\rm{T}}+ V_{\rm{S}} S_{\rm{S}}+ V_B S_{\rm{B}})\}/V_{\rm{IP}},
	\label{eq:totalsalt}
\end{equation}
where $V_{\rm{i}}$ are box volumes. This leaves a 2-dimensional system that models the dynamic change in salinity $S_{\rm{T}}$ of the Tropical (T) and $S_{\rm{N}}$ of the Northern (N) boxes. The strength of the AMOC is linearly proportional to $S_N$ and is given by
$$
q=\frac{\lambda [\alpha (T_{\rm{S}} - T_0) + \beta (S_{\rm{N}} - S_{\rm{S}})]}{1+\lambda\alpha\mu},
$$
where $\lambda, \alpha, \beta$ and $\mu$ are model parameters. The temperature $T_S$ of the Southern box and the global average temperature $T_0$ are assumed to be fixed. 
The system is then given (as in \cite{chapman2025quantifying}) by:
\begin{equation}
\left.
    \begin{array}{rl}
        \frac{dS_{\rm{N}}}{dt}&=f_1(S_{\rm{N}},S_{\rm{T}},\mathcal{H}(t))\\
        \frac{dS_{\rm{T}}}{dt}&=f_2(S_{\rm{N}},S_{\rm{T}},\mathcal{H}(t))
    \end{array}
\right\}
    \label{eq:AMOC-boxmodel}
\end{equation}
where
\begin{equation}
	\label{eq:AOMC_model_q>0_3box}
	\left.\begin{array}{rcl}
		f_1	 &:=& Y[q(S_{\rm{T}} - S_{\rm{N}}) + K_{\rm{N}} (S_{\rm{T}} - S_{\rm{N}}) - F_{\rm{N}}(\mathcal{H}(t)) S_0]/V_{\rm{N}}, \\
		f_2	&:=& Y[q [\gamma S_{\rm{S}} + (1- \gamma) S_{\rm{IP}} - S_{\rm{T}}] + K_{\rm{S}} (S_{\rm{S}} - S_{\rm{T}})\\
		&&+K_{\rm{N}}(S_{\rm{N}} - S_{\rm{T}}) - F_{\rm{T}}(\mathcal{H}(t)) S_0]/V_{\rm{T}},\\
	\end{array}\right\}
\end{equation}
for $q\geq 0$ and
\begin{equation}
	\label{eq:AOMC_model_q<0_3box}
	\left.\begin{array}{rcl}
		f_1 &:=& Y[|q|(S_{\rm{B}} - S_{\rm{N}}) + K_{\rm{N}} (S_{\rm{T}} - S_{\rm{N}}) - F_{\rm{N}}(\mathcal{H}(t)) S_0]/V_{\rm{N}}, \\
		f_2	&:=& Y[|q| (S_{\rm{N}} - S_{\rm{T}}) + K_{\rm{S}} (S_{\rm{S}} - S_{\rm{T}})\\
		&&+K_{\rm{N}}(S_{\rm{N}} - S_{\rm{T}}) - F_{\rm{T}}(\mathcal{H}(t)) S_0]/V_{\rm{T}},\\
	\end{array}\right\}
\end{equation}
for $q<0$.

We refer to \cite{wood2019observable,Alkhayoun2019,chapman2025quantifying} for a discussion of the interpretation and fitting of the parameters to the data from the global circulation model. Specifically, here we use parameter values for doubled $CO_2$, based on FAMOUS$_B$ runs \cite{smith2012famous}, which can be found in Table~\ref{table:1} in Appendix~\ref{app:par_vals}. As in these papers, we assume that freshwater fluxes are linearly proportional to a hosing parameter, $\mathcal{H}$
\begin{equation}
F_{\rm{N}}(\mathcal{H}(t)) = F_{\rm{N,0}} + F_{\rm{N,1}}\mathcal{H}(t), \qquad F_{\rm{T}}(\mathcal{H}(t)) = F_{\rm{T,0}} + F_{\rm{T,1}}\mathcal{H}(t).
\end{equation}
This freshwater hosing parameter acts as the forcing parameter, which if increased too much or too quickly can lead to a collapse of the AMOC that for this model corresponds to a reversal in the flow of water.

We will consider linear piecewise hosing profiles, $\mathcal{H}(t)$ given by
\begin{equation}
\mathcal{H}(t) =  
    \begin{cases} 
      H_0 & t\leq 0, \\
      H_0 + (H_{\max}-H_0)t/T_{rise} & 0\leq t\leq T_{rise}, \\
      H_{\max} & T_{rise}\leq t \leq T_{rise}+T_{plat}, \\
      H_{\max} - (H_{\max}-H_0)(t-T_{rise}-T_{plat})/T_{fall} & T_{rise}+T_{plat} \leq t \leq T_{rise}+T_{plat}+T_{fall}, \\
      H_0 & t\geq T_{rise}+T_{plat}+T_{fall}.
    \end{cases}
\label{eq:hosing_profile}
\end{equation}
The forcing starts at some initial level ${H}_0$ and increases linearly to a maximum level $H_{\max}$ over a duration $T_{rise}$. The forcing plateaus at this maximal level for a duration $T_{plat}$, before linearly decreasing back to $H_0$ over a duration $T_{fall}$ and then remains at $H_0$ thereafter. The default forcing parameter values are provided in Table~\ref{table:2} in Appendix~\ref{app:par_vals}. 

\begin{figure}
\centering
\includegraphics[width=0.5\linewidth]{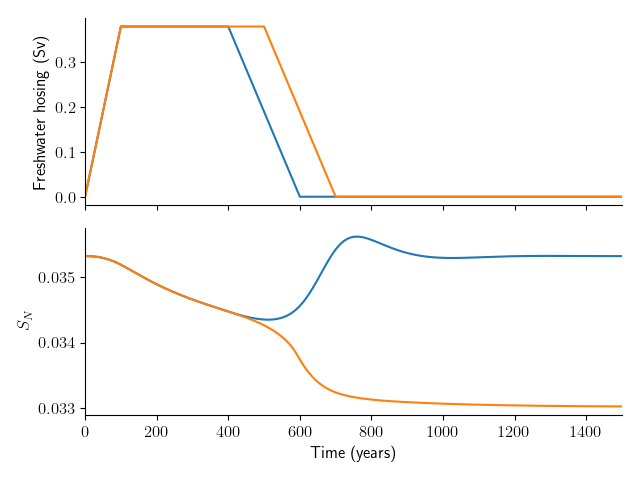}%
\includegraphics[width=0.5\linewidth]{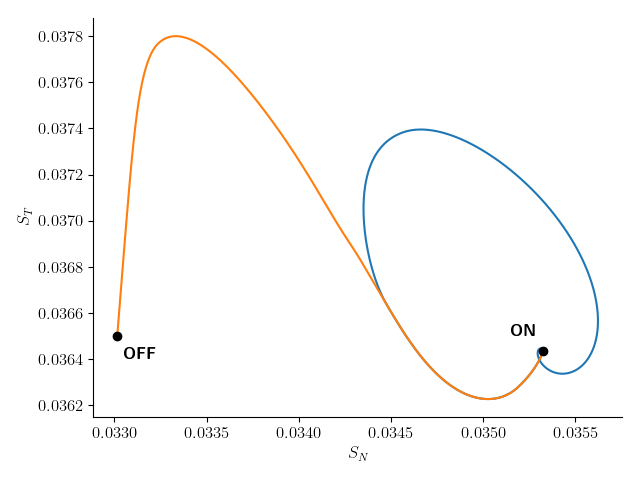}%
\caption{\textbf{Rate-induced tipping in the AMOC 3-box model.} (a) Two time series profiles for freshwater hosing, defined by Eqn.~\eqref{eq:hosing_profile}. The forcing profiles consist of a $T_{rise} = 100$ year linear ramp up to a maximal level $H_{\max}$, which plateaus for $T_{plat} = 300$ years (blue) or 400 years (orange) before a linear ramp down back to zero hosing ($H_0 = 0$) over a $T_{fall} = 200$ year duration. (b) Time series of the salinity of the northern box in response to the freshwater hosing profiles given in (a). (c) Trajectories in the phase plane of the salinities in the tropical and northern boxes given the freshwater hosing forcing in (a). Note that the edge state is off the plot (to the top left).
}
\label{fig:R_tipping_AMOC}
\end{figure}

Importantly, we consider forcing scenarios (with fixed $H_0$ and $H_{\max}$) such that no critical thresholds are crossed, but the system can still undergo R-tipping depending on the durations of $T_{rise}$, $T_{plat}$ and $T_{fall}$. Figure~\ref{fig:R_tipping_AMOC} demonstrates this for two profiles (plotted in panel (a)) with different plateau durations, $T_{plat}$, but with the same $T_{rise}$ and $T_{fall}$. Panel (b) shows that the salinity levels in the Northern box remain high for the blue scenario, corresponding to the AMOC avoiding tipping. However, if the reversal in the hosing is delayed by an extra 100 years, this would be sufficient for the Northern box to become much fresher, signifying an AMOC collapse (orange trajectory). 

\begin{figure}
\centering
\includegraphics[width=0.5\linewidth]{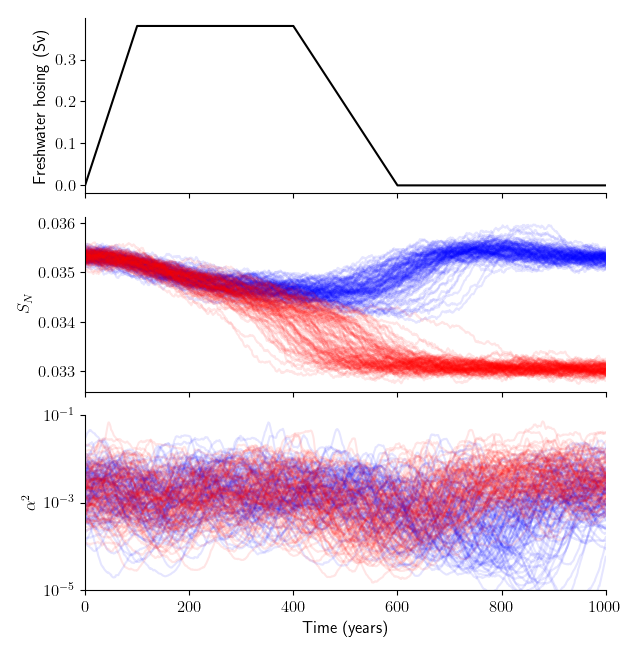}
\caption{\textbf{No early warning via the return rate for ensembles of tipping and non-tipping trajectories.} (a) Time series of a sample freshwater hosing profile \eqref{eq:hosing_profile} with $T_{rise} = 100$, $T_{plat} = 300$, $T_{fall} = 200$, $H_0 = 0$ and $H_{\max} = 0.38$. (b) Monte-Carlo simulations of 100 stochastic trajectories (using \eqref{eq:AMOC-boxmodel_noise} with $\sigma=0.01$, $A_{11} = 0.1263$, $A_{12} = -0.0869$, $A_{21} = 0$, $A_{22} = 0.1008$) that do not tip (blue) and 100 that do tip (red) for the salinity in the northern box. (c) The square of the return rate calculated on (b) on a sliding window of length 200 years.}
\label{fig:R_tipping_return_rate_ensemble}
\end{figure}

The trajectories of these two scenarios have been plotted in the phase plane of the salinities of the Tropical and Northern boxes in panel (c). 
Increasing the freshwater hosing, dramatically reduces the basin of attraction of the ON state and reduces its stability. We first consider whether the return rate (representing traditional EWS) can identify this declining stability. In panel (b) of Figure~\ref{fig:R_tipping_return_rate_ensemble} we plot the salinity in the Northern box for an equal number of stochastic trajectories that tip and do not tip for the forcing profile in panel (a), which does not give tipping without noise. The square of the return rate measured on the salinity in the Northern box is plotted in panel (c). This shows no clear trends of declining resilience, let alone being able to separate the coloured trajectories. Thus, with little evidence of traditional early warning, we will now focus on the prospect of geometric early warning signals, using the signed distance to thresholds in phase space. A particular focus will be on seeing if the R-tipping indicator can offer more insight. 

By definition, a deterministic trajectory will never cross the R-tipping threshold, as shown in the top row in Figure~\ref{fig:R_tipping_edge_state_threshold}. Panel (a) shows the time evolution of the R-tipping threshold for the blue freshwater hosing profile given in Figure~\ref{fig:R_tipping_AMOC}(a). The AMOC ON state, which is where the system is initialised, is within the small region bounded by the R-tipping threshold at the start of the forcing. Therefore, the blue trajectory will not tip and successfully tracks the ON state. However, for the orange forcing given in Figure~\ref{fig:R_tipping_AMOC}(a), the R-tipping threshold encloses a much smaller region at the start of the forcing, which does not include the ON state; see Figure~\ref{fig:R_tipping_edge_state_threshold}(b). Hence, as the system is outside the region enclosed by the R-tipping threshold, it remains outside and tips to the collapsed, OFF state. 

\subsection{R-tipping with stochastic variability}

Let us now consider system~\eqref{eq:AMOC-boxmodel} but with stochastic variability, i.e.:

\begin{equation}
\left.
    \begin{array}{rl}
        \mathrm{d}S_{\rm{N}} &=f_1(S_{\rm{N}},S_{\rm{T}},\mathcal{H}(t))\mathrm{d}t + \sigma(A_{11}\mathrm{d}W_1 + A_{12}\mathrm{d}W_2),\\
        \mathrm{d}S_{\rm{T}} &=f_2(S_{\rm{N}},S_{\rm{T}},\mathcal{H}(t))\mathrm{d}t + \sigma(A_{21}\mathrm{d}W_1 + A_{22}\mathrm{d}W_2),
    \end{array}
\right\}
    \label{eq:AMOC-boxmodel_noise}
\end{equation}
where a matrix $A$ (with elements $A_{ij}$) provides the noise structure, and similar to~\eqref{eq:SDE-nonaut}, $\sigma$ gives the noise strength (default parameter values provided in Table~\ref{table:2} in Appendix~\ref{app:par_vals}) and $\mathrm{d}W_i$ represents independent Wiener white noise processes. 

\begin{figure}
\centering
\includegraphics[width=0.5\linewidth]{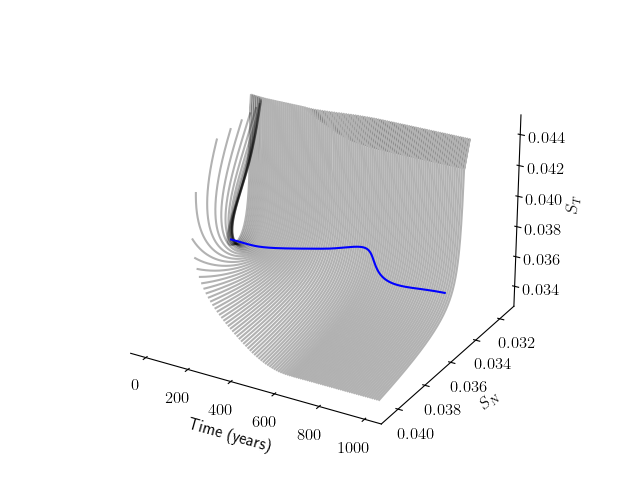}%
\hspace{-4ex}
\includegraphics[width=0.5\linewidth]{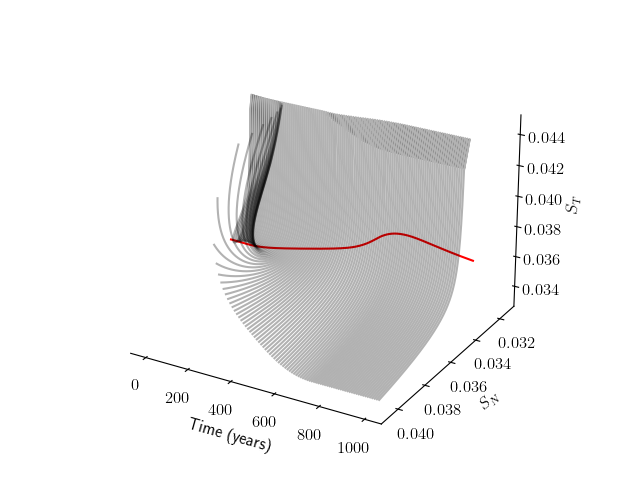}%
\vspace{-0ex}
\includegraphics[width=0.5\linewidth]{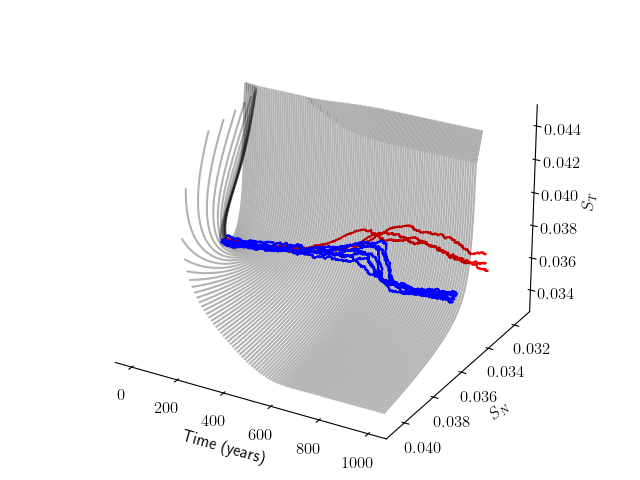}%
\hspace{-4ex}
\includegraphics[width=0.5\linewidth]{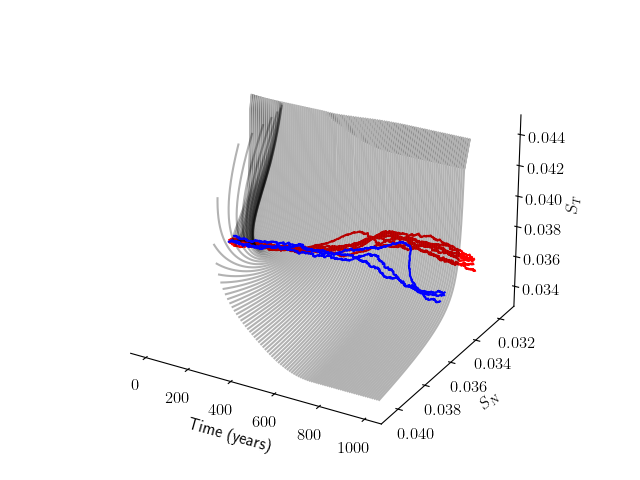}%
\caption{\textbf{Time evolving R-tipping threshold.} Time evolving R-tipping threshold represented as a series of contours in 3D for the piecewise freshwater hosing profiles \eqref{eq:hosing_profile} given in Figure~\ref{fig:R_tipping_AMOC}(a) with $T_{plat} = 300$ years (left column) or $400$ years (right column). Deterministic (top row) and ten sample stochastic (bottom row; using \eqref{eq:AMOC-boxmodel_noise} with $\sigma=0.01$, $A_{11} = 0.1263$, $A_{12} = -0.0869$, $A_{21} = 0$, $A_{22} = 0.1008$) trajectories are plotted in blue for no tipping and red for tipping. Note R-tipping threshold not plotted here prior to the start of forcing, see Figure~\ref{fig:R_tipping_edge_state_threshold_before_forcing} for prior evolution.}
\label{fig:R_tipping_edge_state_threshold}
\end{figure}

The inclusion of stochastic variability means that it is now possible to cross surfaces that would be R-tipping thresholds for the deterministic system. For instance, in Figure~\ref{fig:R_tipping_edge_state_threshold}(c) all trajectories are initialised at the ON state, which is inside the region bounded by the R-tipping threshold prior to the start of the forcing. However, some trajectories (coloured red) tip to the collapsed state. Similarly, in Figure~\ref{fig:R_tipping_edge_state_threshold}(d) all trajectories are initialised outside the region bounded by the R-tipping threshold, but there are some blue trajectories that avoid tipping, ending inside the region bounded by the R-tipping threshold. However, the further the system travels beyond the R-tipping threshold, the more likely the system will tip.

\begin{figure}
\centering
\includegraphics[width=0.5\linewidth]{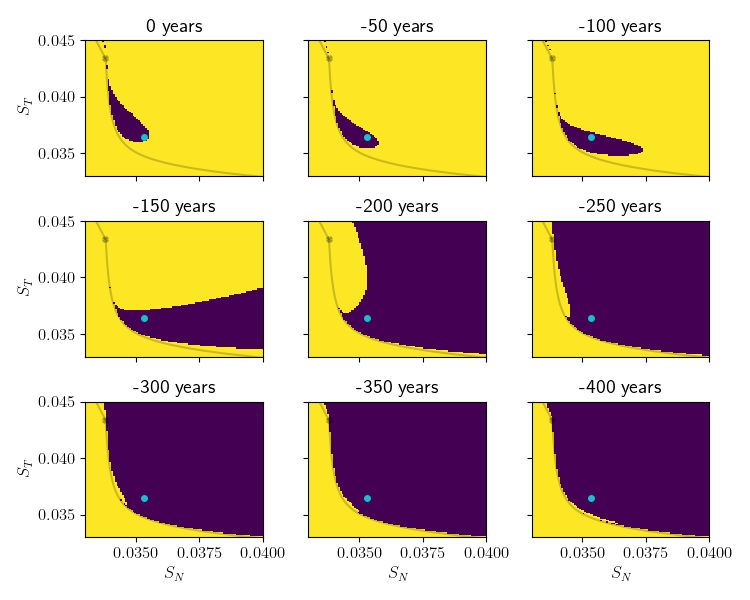}%
\includegraphics[width=0.5\linewidth]{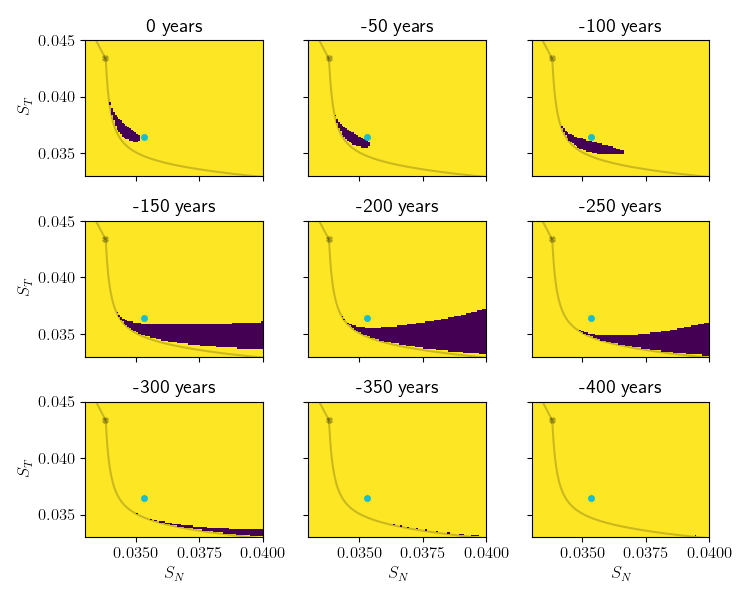}%
\caption{\textbf{Points in phase space that tip based on initialisation period prior to start of forcing.} At the number of years prior to the forcing starting if the system is at a yellow point it tips to the OFF state or a purple point it returns to the ON state (light blue dot) for the piecewise freshwater hosing profiles given in Figure~\ref{fig:R_tipping_AMOC}(a) with $T_{plat} = 300$ years (left) or $400$ years (right). The edge state at no freshwater hosing is given by a grey cross and its stable manifolds by the grey curves.}
\label{fig:R_tipping_grid_before_forcing}
\end{figure}

It is important to note that, unlike the basin boundary, the R-tipping threshold will change in the lead up to the starting of the forcing. This can first be established from Figure~\ref{fig:R_tipping_grid_before_forcing}, which depicts regions of the phase plane that deterministically avoid tipping (purple points) and regions of tipping (yellow points) for initialising the system at these points the given number of years before the forcing starts. The left set of panels corresponds to a forcing plateau of $T_{plat} = 300$ years, and there is only a small region that avoids tipping. Crucially, this region does include the stable ON state, and thus, for longer initialisation times before the start of the forcing, this purple region expands.  For a sufficiently long time prior to the forcing starting, this region will converge to the basin of attraction for the ON state. This is because, given sufficient time, all points will converge to the ON state by the time the forcing starts, and starting at the ON state avoids tipping for this forcing profile.

The right-hand set of panels are for a forcing plateau of $T_{plat} = 400$ years, where starting at the ON state when the forcing starts causes tipping. There is, however, a small region that, if the system was located at, would avoid tipping to the OFF state. 
Upon increasing the number of years before starting the forcing, this region initially grows. For example, 200 years before the forcing starts, there is a relatively large area of the phase plane that will converge to the smaller region when the forcing starts. However, for a sufficiently early initialisation, the region shrinks and disappears (all points are yellow), as all points within the basin of attraction will converge to the ON state, prior to the forcing starting. For this particular forcing profile, starting at the ON state will result in tipping, and therefore no starting position in the phase plane (apart from very close to the basin boundary) can avoid tipping. 

At the critical plateau duration for starting at the ON state, approximately periodic behaviour occurs in the lead up to the start of the forcing, with points in phase space frequently switching between tipping and not tipping; see animation provided in Supplementary Material. Figure~\ref{fig:R_tipping_edge_state_threshold_before_forcing} provides some snapshots of the evolution of the R-tipping threshold, which should encompass the purple points from Figure~\ref{fig:R_tipping_grid_before_forcing}, prior to the start of the forcing. 
Here it is clearer to see that although the region shrinks as time increases (time before forcing decreases), the region also shifts location. Therefore, even if the system is initialised in the upper segment of the region enclosed by the blue curve (for the forcing corresponding to the left panel), but 150 years prior to the forcing starting, tipping would still occur. This is because, without forcing, the system spirals to the ON state; however, when the forcing starts, the rotation may mean that the system is no longer within the blue region, despite converging to the ON state. 

\begin{figure}
\centering
\includegraphics[width=0.5\linewidth]{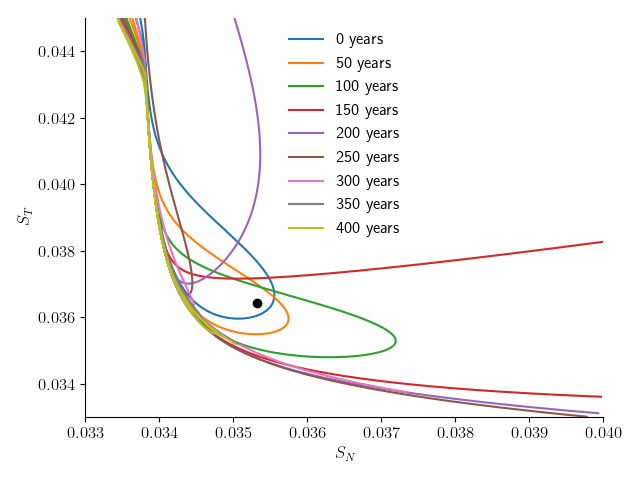}%
\includegraphics[width=0.5\linewidth]{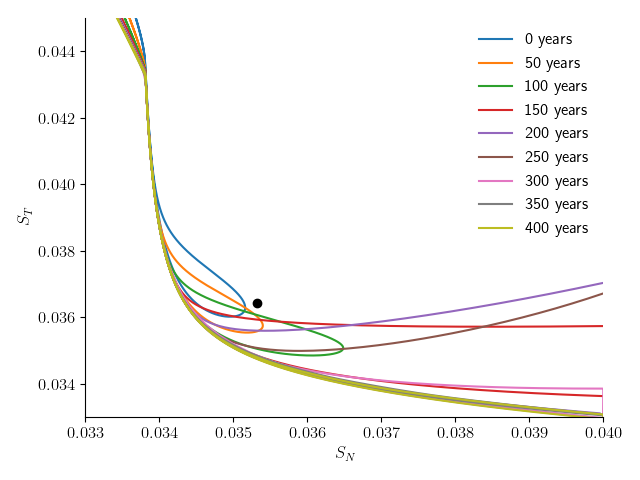}%
\caption{\textbf{Time evolution of R-tipping threshold prior to the start of the forcing.} R-tipping threshold plotted at different time intervals prior to the start of the piecewise freshwater hosing profiles given in Figure~\ref{fig:R_tipping_AMOC}(a) with $T_{plat} = 300$ years (left panel) or $400$ years (right panel). The stable ON state is given by the black dot.}
\label{fig:R_tipping_edge_state_threshold_before_forcing}
\end{figure}

For a longer duration of peak freshwater hosing ($T_{plat}$), the ON state is outside the R-tipping threshold region at the start of the forcing; see Figure~\ref{fig:R_tipping_edge_state_threshold_before_forcing}(b). Therefore, the ON state cannot be within the R-tipping threshold at any stage (as apparent from Fig~\ref{fig:R_tipping_grid_before_forcing} (b)), and the R-tipping threshold converges to the basin boundary for longer initialisations prior to the start of the forcing.

\section{Quantifying the skill of early warning classifiers}
\label{sec:EWS_skill}

For the 3-box AMOC system with noise \eqref{eq:AMOC-boxmodel_noise}, one can argue that a simpler metric based on the level of salinities in either the Tropical ($S_T$) or Northern ($S_N$) box could diagnose tipping (a similar approach has been taken in \cite{jackson2023understanding,Romanou2023}). Indeed, looking at Figure~\ref{fig:R_tipping_return_rate_ensemble}(b) suggests this is possible, by devising geometrical early warning signals based on the system salinities in either of the two boxes relative to some pre-determined \textit{critical salinity levels}, presumably with reference to the salinities of the unstable edge state. Such salinity thresholds can be used to separate tipping and non-tipping trajectories for an unforced system. 

In the presence of strong forcing, we will need to define time-dependent thresholds that asymptote to the fixed thresholds as the trajectories choose their fate (after forcing is off). We call these the ``optimum thresholds''. Note that the R-tipping threshold itself is an object in the phase space relative to which some fixed and optimum threshold distances can be defined. Keeping this in mind, we compare the three metrics, namely 1) R-tipping indicator, 2) deviation from a critical salinity level in the Northern box, and 3) deviation from a critical salinity in the Tropical box; as geometrical early warning signals defined in terms of ``fixed'' or ``optimum'' thresholds.

\begin{figure}
\centering
\includegraphics[width=0.5\linewidth]{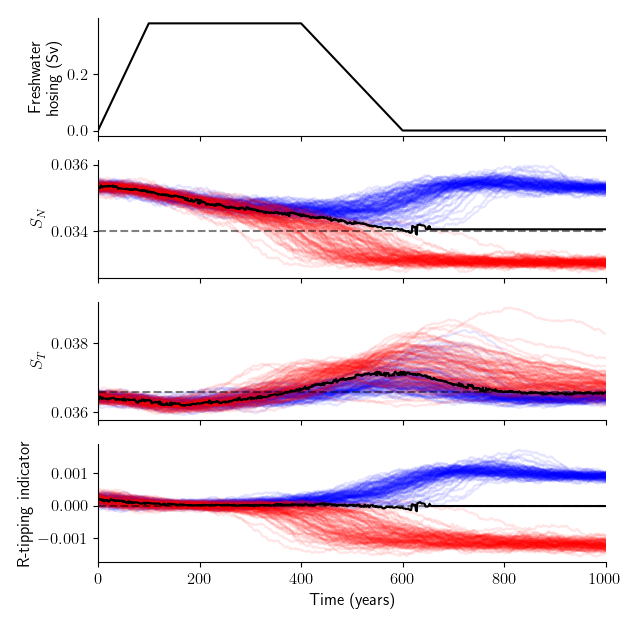}
\caption{\textbf{Geometric indicators for ensembles of tipping and non-tipping trajectories.} (a) Time series of a sample freshwater hosing profile \eqref{eq:hosing_profile} with $T_{rise} = 100$, $T_{plat} = 300$, $T_{fall} = 200$, $H_0 = 0$ and $H_{\max} = 0.38$. 100 sample stochastic trajectories (using \eqref{eq:AMOC-boxmodel_noise} with $\sigma=0.01$, $A_{11} = 0.1263$, $A_{12} = -0.0869$, $A_{21} = 0$, $A_{22} = 0.1008$) that do not tip (blue) and 100 that do tip (red) plotted as time series for the salinity in the northern box (b), salinity in the Tropical box (c) and the R-tipping indicator (d). Black solid curve shows optimal threshold as a function of time and grey dashed line a fixed threshold used for separating tipped and non-tipped trajectories.}
\label{fig:R_tipping_edge_state_threshold_ensemble}
\end{figure}

Figure~\ref{fig:R_tipping_edge_state_threshold_ensemble} plots the time series of the salinities in the Northern and Tropical boxes and the R-tipping indicator threshold for an equal number of sample stochastic trajectories that tip and do not tip for the same forcing profile. Panel (b) shows that once the salinity in the Northern box ($S_{N}$) gets below a certain level (indicated by the grey dashed line), the trajectory always ends up tipping to the collapsed state. However, arguably, this warning would come too late, and so one might ask if there is a more optimal level that can be chosen. The salinity in the Tropical box ($S_{T}$) appears to be even harder to define a threshold for (see panel (c)), as the ON and OFF states have very similar salinity levels as shown by Figure~\ref{fig:R_tipping_AMOC}(c). In the following, we will consider which of these metrics, including the R-tipping indicator (shown in Figure~\ref{fig:R_tipping_edge_state_threshold_ensemble}(d)), is best for providing early warning.
Further, we will determine an optimal threshold that minimises the number of false positives and false negatives.

\subsection{ROC analysis of performance}

Receiver Operating Characteristic (ROC) curves provide a graphical illustration of the performance of an indicator to identify skill of predicted events \cite{fawcett2006introduction}. They classify trajectories in a binary system (in our case, tipping vs non-tipping) according to a threshold, where the optimal threshold is chosen to balance the probability of false positives (horizontal) against the probability of true positives (vertical). 

\begin{figure}
\centering
\includegraphics[width=0.7\linewidth]{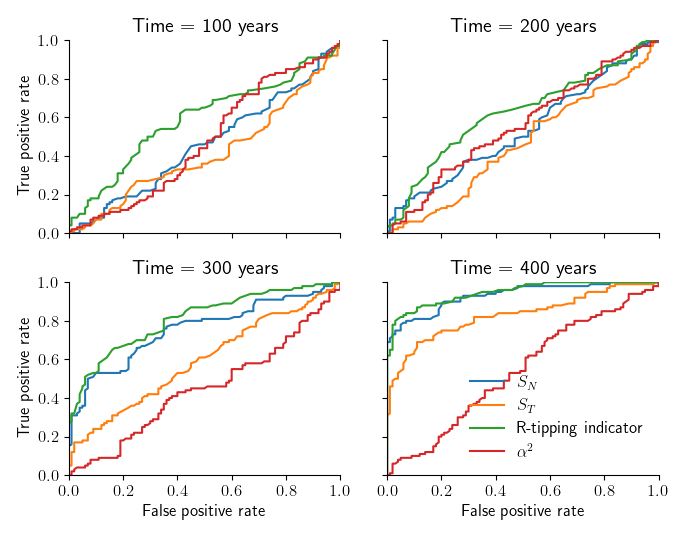}
\caption{\textbf{ROC curves for early detection of tipping.} ROC curves at different time intervals showing the performance of using the Northern box salinity (blue), Tropical box salinity (orange), R-tipping indicator (green), and return rate (red) as indicators for detecting tipping.}
\label{fig:ROC_curves}
\end{figure}

Choosing a threshold too extreme in one direction will result in failure to predict any tipping (zero true positives and zero false positives). As the threshold is changed, the number of false positives and true positives will both monotonically increase until the threshold is such that tipping will always be predicted, resulting in many false positives. The closer the curve approaches $(0,1)$, the better the performance of the indicator. Figure~\ref{fig:ROC_curves} provides ROC curves at different time intervals using the salinity levels in the Northern and Tropical boxes and the R-tipping indicator, as well as the square of the return rate ($\alpha^2$). After the ramp up has finished (100 years, panel (a)) all the ROC curves lie roughly on the 1-1 line indicating the predictability of a coin flip. As time increases, the salinity in the Northern box and R-tipping indicator gain skill (ROC curves move towards (0,1)), with the R-tipping threshold slightly outperforming the Northern box salinity. Notice that this is still 100 years before freshwater starts to be removed (panel (c)), and the salinity in the Tropical box and the return rate still offer very little skill. At the time freshwater is starting to be removed (400 years, panel (d)), the salinity in the Northern box and the R-tipping threshold have converged and provide high skill in offering early warning (ROC curves converge close to (0,1)). The salinity in the Tropical box also now offers a reasonable amount of skill, yet the return rate is not informative. 


\begin{figure}
\centering
\includegraphics[width=0.5\linewidth]{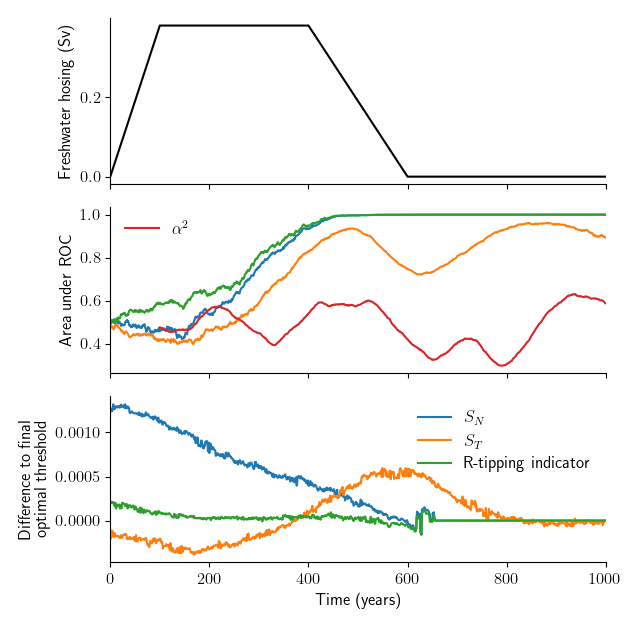}
\caption{\textbf{Time evolution for the area under ROC curves for early detection of tipping.} Time series of (a) freshwater hosing forcing, (b) area under the ROC curve, (c) change in optimal threshold for salinity in the Northern box (blue), salinity in the Tropical box (orange), R-tipping indicator (green), and return rate (red). Note the optimal threshold for the return rate is not plotted as it has different units to the other three thresholds.}
\label{fig:AUC}
\end{figure}

This can be captured more accurately by measuring the area under the ROC curve (AUC) \cite{fawcett2006introduction}. If the AUC is 0.5, then the predictor is only as good as a coin flip, while an AUC = 1 reflects a perfect indicator. Figure~\ref{fig:AUC}(b) plots a time series of the AUC for the four different indicators. As already identified from the ROC curves, while the forcing is linearly increasing, the AUC for all indicators is around 0.5. Despite the AUC for R-tipping threshold increasing marginally earlier than the AUC for the salinity in the Northern box there is little difference between their overall skill. By the start of the freshwater hosing decreasing, the AUC's have converged close to one signifying a clear separation between tipping and non-tipping trajectories in either the salinity of the northern box or the R-tipping indicator. The increase in AUC for the salinity in the Tropical box is delayed by around 100 years but then increases at a similar rate. However, the AUC or skill peaks during the middle of the decline in freshwater
(see Figures~\ref{fig:R_tipping_AMOC}(c) and \ref{fig:R_tipping_edge_state_threshold_ensemble}(c)) before relaxing back to similar levels of non tipping trajectories, which causes a drop in skill of the indicator. The AUC for the return rate always hovers around 0.5, but then even drops to being worse than a coin flip when the forcing finishes. 

\subsection{Optimal thresholds for early warnings}

An optimal threshold  can be identified as the threshold that simultaneously maximises the number of true positives and minimises the number of false positives. A rational choice of this threshold is to choose the shortest distance from the ROC curve to $(0,1)$, though there are clearly other ways to do this depending on how false positives and false negatives are weighted \cite{fawcett2006introduction}. The optimal threshold will change over time and is plotted as a black curve in Figure~\ref{fig:R_tipping_edge_state_threshold_ensemble}. Noticeably, the optimal threshold for R-tipping indicator changes much less than those for either salinity in the Northern or Tropical boxes, as also shown by Figure~\ref{fig:AUC}(c). This highlights the additional benefit of using the R-tipping threshold, as in practice a fixed threshold would be chosen. 

\begin{figure}
\centering
\includegraphics[width=0.5\linewidth]{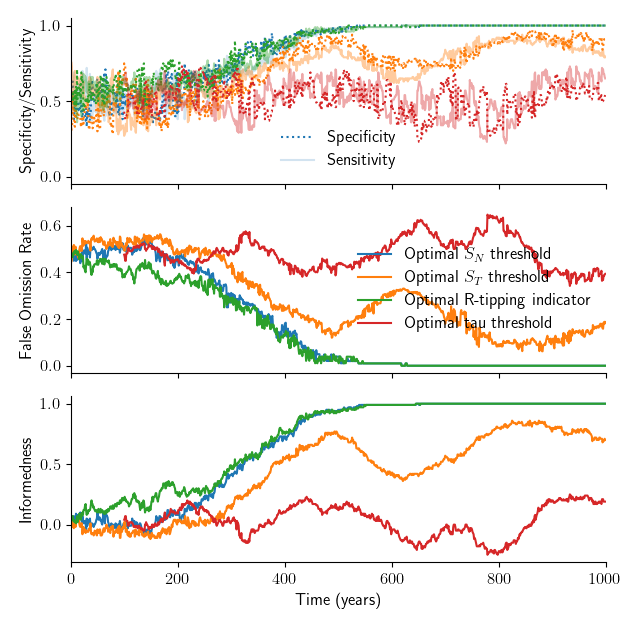}%
\includegraphics[width=0.5\linewidth]{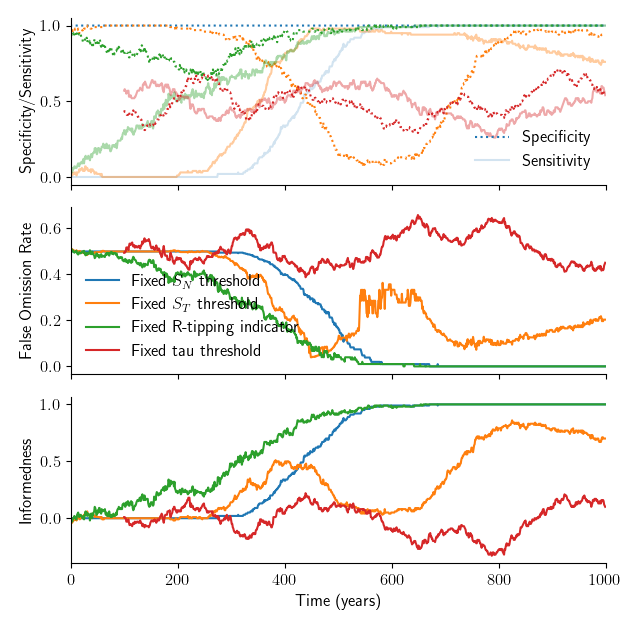}%
\caption{\textbf{Performance statistics of using the optimal and a fixed threshold in both the northern box salinity and signed distance to R-tipping threshold.} Time series of the sensitivity (TPR; solid curve) and specificity (1-FPR; dotted curve) in the top row; false omission rate in the middle row; Informedness (Sensitivity + Specificity - 1) in the bottom row for the optimal threshold in the left column and for a fixed threshold in the right column using the salinity in the Northern box (blue), salinity in the Tropical box (orange), R-tipping indicator (green), and return rate (red).}
\label{fig:indicator_stats}
\end{figure}

Figure~\ref{fig:indicator_stats} plots various statistics that measure the performance of each indicator both for the optimal threshold (left column) and a fixed threshold (right column). The fixed threshold (shown as a grey dashed line in Figure~\ref{fig:R_tipping_edge_state_threshold_ensemble} for the geometric indicators) is chosen such that at the earliest possible stage there is a perfect classification between tipped and non-tipped trajectories. For the R-tipping indicator, a natural and sensible choice is a distance of zero, whereas for the salinities in the two boxes it is a bit more arbitrary (and might change for different models/parametrizations), chosen here to be $0.034$ and $0.0366$ for the salinities in the Northern and Tropical boxes, respectively.

Given the ROC curves it is unsurprising that the salinity in the Northern box and the R-tipping indicator perform equally well if the optimal threshold is known. However, the advantage of the R-tipping indicator can be observed by comparing the fixed threshold statistics. For the R-tipping indicator, the statistics for a fixed threshold of zero are very similar to using the optimal threshold. On the other hand, the salinity in the Northern box performs much worse with a fixed threshold than with the use of the optimal threshold. Up until the forcing beginning to decline the salinity in the Tropical box even outperforms using the salinity in the Northern box.

\begin{figure}
\centering
\includegraphics[width=0.5\linewidth]{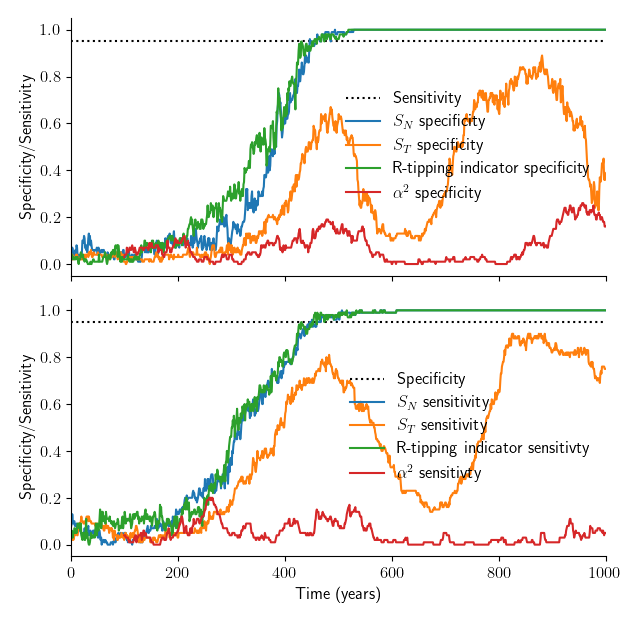}
\caption{\textbf{Sensitivity and specificity performance with the other metric set to a minimum value.} (a) Time series of the specificity (coloured curves) given at least 0.95 sensitivity (black dotted line). (b) Time series of the sensitivity (coloured curves) given at least 0.95 specificity (black dotted line). Statistics calculated based on the salinity in the Northern box (blue), salinity in the Tropical box (orange), R-tipping indicator (green), and return rate (red).}
\label{fig:max_FP_FN}
\end{figure}

Extending this, one might be interested in limiting either the number of false positives or false negatives, as performed in Figure~\ref{fig:max_FP_FN}. Panel (a) displays the specificity for the four metrics while choosing a threshold that ensures a minimum sensitivity of 0.95 (i.e. a maximum 5\% false negative rate). The R-tipping indicator marginally outperforms the salinity in the northern box for the first 400 years, before both converge to a sensitivity of one. However, the rate of false positives for using the salinity in the Tropical box fluctuates over time. The square of the return rate though performs even worse with the specificity rarely above 0.2. Near identical behaviour is found in all four metrics if we were to instead limit the number of false positives; see panel (b).  

\section{Discussion}
\label{sec:Discussion}

Critical slowing down and classical early warning signals are of limited use for understanding tipping in rapidly forced systems. We therefore propose that the signed distance of the system state to its R-tipping threshold (a geometrical measure we call the R-tipping indicator) can provide early warning of tipping in this case. The R-tipping threshold deterministically separates trajectories in phase space that undergo tipping to different attractors, based on the future profile of the forcing. This method requires knowledge of the system's dynamics, its current location in phase space, and the future profile of the forcing. However, no knowledge of the system's past is required, specifically how it reached its current state and the prior forcing. 

We demonstrate the applicability of such a geometric early warning signal for a box model for the Atlantic Meridional Overturning Circulation (AMOC), extending the suggestion in \cite{chapman2025quantifying} to estimate this eventual fate even while the forcing is active. Knowledge of the future forcing, system dynamics and edge state allows us to distinguish between non-tipping and tipping trajectories. We verify that in the presence of stochastic terms, this prediction which gains skill as time progresses. It will be interesting to compare this still to machine learning approaches such as those developed in \cite{huang2024anticipating}.

Gaining a skilful early warning, even before freshwater hosing (forcing) starts to decline, has implications for mitigation. For example, actions can be taken to increase mitigation efforts if it is found that the proposed rate of reversal of forcing would not be sufficient to prevent tipping. 

Alternatively, an earlier warning will offer the opportunity to implement adaptation methods. Such methods often rely on knowledge of the future fate of the system. Consequently, the ultimate goal of early warning signals should be to determine whether a system will tip in the future or not. However, classical methods such as those used to detect critical slowing down by design can only in the best case inform if a system is approaching a tipping point and cannot account for the future forcing, and hence the future fate of the system. In such cases, geometric early warning signals that take into account the system dynamics, its current state, and future forcing offer a possible solution. 

This could be a simple metric that defines a critical threshold in a system variable and issues a warning based on whether the threshold is crossed \cite{wood2019observable}. In the AMOC example, we find that salinity in the Northern box offers a skilful early warning. However, it is not always apparent which variable or observable has the most potential for skilful early warning. \cite{lohmann2025role} found that looking for critical slowing down in the direction of the edge state can offer greatest skill. 

The R-tipping indicator can be viewed as 
an optimal (nonlinear) combination of the state variables.
More precisely, the R-tipping indicator is based on which side the system lies relative to the threshold. We provide one way to identify these optimum thresholds by analysing the ROC curves. This analysis captures the skill of the proposed geometric EWS and can be used to identify optimum thresholds that maximize \textit{informedness}.

There could be limitations on knowledge of future forcing in most realistic situations. Different types of time-varying forcing may exist \cite{ashwin2025early}, or there may be uncertainties in the assumed linear profile. Even though the proposed geometric early warning signal offers a skilful warning, uncertainties in the future forcing will ultimately lead to uncertainties in the tipping behaviour. So, a framework that allows for probabilities of tipping to account for such uncertainties needs to be developed. 

Diagnosing the R-tipping indicator for complex models or with added complexity of the edge state presents a particular limitation of this approach. Data-driven model discovery utilising deep learning methods might help overcome some of these limitations. As more observational datasets are acquired, deep learning and artificial intelligence systems can enhance mathematical modelling. For existing high-dimensional complex system models, such as those in climate, neuroscience, or infrastructure (e.g. power grids), the limitations are mostly in computation. This remains an exciting open area of research in computer science and complex systems. Geometric early warning signals such as the proposed R-tipping indicator, offer a quantitative framework towards designing a realistic early warning system.

\section*{Data Availability Statement}

The code to reproduce the results presented in this paper is available via \url{https://github.com/PaulRitchie89/R-tipping_threshold}.

\section*{Acknowledgements}

PR, SK, and PA thank AdvanTip and PR and PA thank ClimTip for funding support. The ClimTip project has received funding from the European Union's Horizon Europe research and innovation programme under grant agreement No. 101137601: Funded by the European Union. Views and opinions expressed are however those of the author(s) only and do not necessarily reflect those of the European Union or the European Climate, Infrastructure and Environment Executive Agency (CINEA). Neither the European Union nor the granting authority can be held responsible for them. The Advanced Research and Invention Agency (ARIA) fund the AdvanTip project (grant no. SCOP-PR01-P003).

For the purpose of open access, we have applied a Creative Commons Attribution (CC BY) license to any Author Accepted Manuscript version arising from this submission.

\newpage

\bibliographystyle{plain}
\bibliography{refs}

\newpage
\appendix

\section{Numerical approximation of the R-tipping threshold}
\label{app:R-threshold_calc}

We briefly outline the process for calculating the time evolution of the R-tipping threshold for the AMOC model used in this paper; see the Data Availability Statement for access to the code. At the end of the forcing, the R-tipping threshold, by definition, is the boundary that separates the basins of attraction for the ON and OFF states (i.e. the stable manifold of the edge state). The R-tipping threshold is then calculated by numerically integrating backward in time points evenly interpolated along this curve. However, many of these points will rapidly go outside the region in the phase space of interest. Therefore, at small time intervals, the R-tipping threshold is again re-interpolated to generate a new set of points evenly distributed along the curve. This process is repeated until the start of the simulation is reached. It is important to note that the R-tipping threshold will still change in the time leading, up to any changes of the external forcing, as shown in Figure~\ref{fig:R_tipping_edge_state_threshold_before_forcing}.

To estimate a R-tipping indicator, we need to determine algorithmically 
whether the system is inside the region of phase space bounded by the R-tipping threshold.
At any given time, we can determine where the system is in relation to the R-tipping threshold by constructing a linear curve that connects the tipped state to the current location of the system in phase space. If the number of intersections of this curve with the R-tipping threshold is odd, then the system is inside the region bounded by the R-tipping threshold and has a positive signed distance; otherwise, it is outside and has a negative signed distance. 

\section{Parameter values for AMOC model}
\label{app:par_vals}

For computational reasons, the model \eqref{eq:AMOC-boxmodel} (and the stochastic version \eqref{eq:AMOC-boxmodel_noise}) is solved for rescaled salinities, $\tilde{S}_i$:
\begin{equation}
    \tilde{S}_i = 100(S_i-S_0), \qquad i \in \{N, T, S, IP, B\},
\end{equation}
where $S_0$ is a reference salinity. The rescaled salinity for the Indo-Pacific box using \eqref{eq:totalsalt} is therefore given by
\begin{equation}
    \tilde{S}_{IP} = 100(C - (V_N\tilde{S}_N+V_T\tilde{S}_T+V_S\tilde{S}_S+V_B\tilde{S}_B)/100 - S_0(V_N+V_T+V_S+V_{IP}+V_B))/V_{IP}.
\end{equation}
In the main paper, the default model parameter values used are given in Table~\ref{table:1} and the forcing and noise parameter values in Table~\ref{table:2}. 

\begin{table}[h!]
\centering
\begin{tabular}{c | p{9cm}| c | c} 
Parameter & Interpretation & Value & Units \\ 
\hline
$S_0$ & Reference salinity & 0.035 & - \\
$\alpha$ & Temperature coefficient for AMOC strength & 0.12 & $kg\, m^{-3}\, ^{\circ}C^{-1}$  \\
$\beta$ & Salinity coefficient for AMOC strength & 790 & $kg\, m^{-3}$ \\
$\lambda$ & Relationship between density and AMOC strength & $1.62\times 10^7$ & $m^6\,kg^{-1}\,s^{-1}$ \\
$\mu$ & Relationship between Northern box temp and AMOC strength & $22\times 10^{-8}$ & $^{\circ}C\, s\, m^{-3}$ \\
$T_s$ & Southern box temperature & 7.919 & $^{\circ}C$ \\
$T_0$ & Global average temperature & 3.87 & $^{\circ}C$ \\
$\gamma$ & Proportion of water that follows cold water path & 0.36 & - \\
$Y$ & Seconds in a year & $3.15\times 10^7$ & $s\,yr^{-1}$ \\
$K_N$ & Wind-driven fluxes between Northern and Tropical boxes & $1.762\times 10^6$ & $m^3\,s^{-1}$ \\
$K_S$ & Wind-driven fluxes between Southern and Tropical boxes & $1.872\times 10^6$ & $m^3\,s^{-1}$ \\ 
$S_S$ & Southern box salinity & 0.034427 & - \\
$S_B$ & Bottom box salinity & 0.034538 & - \\
$C$ & Total salinity content & $4.4735\times 10^{16}$ & $m^3$ \\
$V_N$ & Northern box volume & $0.3683\times 10^{17}$ & $m^3$ \\
$V_T$ & Tropical box volume & $0.5418\times 10^{17}$ & $m^3$ \\
$V_S$ & Southern box volume & $0.6097\times 10^{17}$ & $m^3$ \\
$V_{IP}$ & Indo-Pacific box volume & $1.4860\times 10^{17}$ & $m^3$ \\
$V_B$ & Bottom box volume & $9.9250\times 10^{17}$ & $m^3$ \\
$F_{N,0}$ & Northern box freshwater constant coefficient & $0.4860\times 10^6$ & $m^3\,s^{-1}$ \\
$F_{N,1}$ & Northern box freshwater linear coefficient & $0.1311\times 10^6$ & $m^3\,s^{-1}$ \\
$F_{T,0}$ & Tropical box freshwater constant coefficient & $-0.997\times 10^6$ & $m^3\,s^{-1}$ \\
$F_{T,1}$ & Tropical box freshwater linear coefficient & $0.6961\times 10^6$ & $m^3\,s^{-1}$ \\
\end{tabular}
\caption{List of default model parameters used for the AMOC model. Values obtained from \cite{Alkhayoun2019} correspond to doubled $CO_2$, based on FAMOUS$_B$ runs \cite{smith2012famous}}
\label{table:1}
\end{table}

\begin{table}[h!]
\centering
\begin{tabular}{c | c| c | c} 
Parameter & Meaning & Value & Units \\ 
\hline\hline
$H_0$ & Initial freshwater hosing level & 0 & - \\
$H_{\max}$ & Maximum freshwater hosing level & 0.38 & - \\
$T_{rise}$ & Linear ramp up duration of freshwater hosing & 100 & $yr$ \\
$T_{fall}$ & Linear ramp down duration of freshwater hosing & 200 & $yr$ \\
$T_{plat}$ & Duration of freshwater hosing at maximum level & 300 & $yr$ \\
\hline
$\sigma$ & Noise strength & 0.01 &  - \\
$A_{11}$ & Local noise variability on Northern box & 0.1263 & $yr^{-1/2}$ \\
$A_{12}$ & Remote noise variability noise on Northern box & -0.0869 &  $yr^{-1/2}$ \\
$A_{21}$ & Remote noise variability on Tropical box & 0 & $yr^{-1/2}$ \\
$A_{22}$ & Local noise variability on Tropical box & 0.1008 & $yr^{-1/2}$ \\
\hline\hline
\end{tabular}
\caption{List of default freshwater forcing and noise parameters used for the AMOC model. Noise parameter values obtained from \cite{chapman2025quantifying} correspond to the 3-box model calibrated against a simulation of the HadGEM3MM model \cite{Chapman2024_noise_est}.}
\label{table:2}
\end{table}

\end{document}